\documentclass[a4j,12pt,reqno]{amsart}

\usepackage{amscd,epsfig,amsthm,graphicx}
\usepackage{amssymb}
\usepackage{amsmath}
\usepackage{amsthm}
\usepackage{amsfonts}
\usepackage[english]{babel}
\usepackage{epic,eepic}
\usepackage[flushleft]{paralist} 
\usepackage[utf8]{inputenc}
\usepackage{enumerate}
\usepackage{graphicx}
\usepackage[all]{xy} 

\makeindex			 
\newtheorem{theorem}{Theorem}

\theoremstyle{definition}
			{Definition}

\newtheorem{remark}{Remark}

\makeatletter
    
    \@addtoreset{equation}{section}
\makeatother

\def\N{{\mathbb N}}
 
 \def\R{{\mathbb R}}

\makeatletter
\def\@setcopyright{}
\def\serieslogo@{}
\makeatother

\begin{document}

%


\author{Takashi MIYAGAWA}      %


\title{Approximate functional equations for the Hurwitz and Lerch zeta-functions}

\begin{abstract}
 As one of the asymptotic formulas for the zeta-function, Hardy and Littlewood gave asymptotic
 formulas called the approximate functional equation.
 In 2003, R. Garunk\v{s}tis, A. Laurin\v{c}ikas, and J. Steuding (in \cite{GLS}) proved 
 the Riemann-Siegel type of the approximate functional equation for the Lerch zeta-function 
 $ \zeta_L (s, \alpha, \lambda ) = \sum_{n=0}^\infty e^{2\pi i n \lambda}(n + \alpha)^{-s} $.
 In this paper, we prove another type of approximate functional equations for the Hurwitz and
 Lerch zeta-functions.
 R. Garunk\v{s}tis, A. Laurin\v{c}ikas, and J. Steuding (in \cite{GLS2}) obtained the
 results on the mean square values of $ \zeta_L (\sigma + it, \alpha , \lambda) $ with respect to $ t $.
 We obtain the main term of the mean square values of $ \zeta_L (1/2 + it, \alpha , \lambda) $ 
 using a simpler method than their method in \cite{GLS2}.
\end{abstract}


\maketitle

\dedicatory{}

\date{}


\maketitle


\section{Introduction and the statement of results}
Let $ s = \sigma + it $ be a complex variable, and let 
$ 0 < \alpha \leq 1, 0 < \lambda \leq 1 $ be real parameters. 
The Hurwitz zeta-function $ \zeta_H(s, \alpha) $ and the Lerch zeta-function
$ \zeta_L(s, \alpha, \lambda) $ are defined by
\begin{equation}\label{Hurwitz-zeta}
		\zeta_H (s, \alpha) = \sum_{n=0}^\infty \frac{1}{(n +\alpha)^s},
\end{equation}
\begin{equation}\label{Lerch-zeta}
		\zeta_L (s, \alpha , \lambda)
               = \sum_{n = 0}^\infty
                 \frac{e^{2 \pi i n \lambda}}{(n + \alpha)^s},
\end{equation}
respectively. There series are absolutey convergent for $ \sigma > 1 $.
Also, if $ 0 < \lambda < 1 $, then the series (\ref{Lerch-zeta}) is convergent even for
$ \sigma > 0 $.

As a classical asymptotic formula for the Riemann zeta-function,
the following was proved by Hardy and Littlewood ({\S}4 in \cite{Tit});
we suppose that $ \sigma_0 > 0 $, $ x \geq 1 $, then
\[
		\zeta(s) = \sum_{n \leq x} \frac{1}{n^s} - \frac{x^{1-s}}{1-s} 
		+ O(x^{-\sigma})
\]	
uniformly for $ \sigma \geq \sigma_0 ,\ |t| < 2\pi x/C $, where $ C > 1 $ 
is a constant. 
Also, Hardy and Littlewood proved the following asymptotic formula 
({\S}4 in \cite{Tit}); we suppose that 
$ 0 \leq \sigma \leq 1,\ x \geq 1,\ y\geq 1 $ and $ 2\pi xy = |t| $, then
\begin{equation}\label{AFE}
	\zeta(s) = \sum_{n \leq x} \frac{1}{n^s} + \chi(s) \sum_{n \leq y} \frac{1}{n^{1-s}}
	+ O(x^{-\sigma}) + O(|t|^{1/2-\sigma}y^{\sigma-1}),
\end{equation}
where $ \chi(s) = 2 \Gamma(1-s) \sin{(\pi s/2)}(2\pi)^{s-1} $ and note that
$ \zeta(s) = \chi(s) \zeta(1-s) $ holds.
This is called approximate functional equation.

Further, there is a Riemann-Siegel type of the approximate functional equation for $ \zeta(s) $;
suppose that $ 0 \leq \sigma \leq 1 , x = \sqrt{t/2\pi} $, and $ N < Ct $ with a sufficiently 
small constant $ C $. Then
\begin{eqnarray}
	&& \zeta(s) = \sum_{n \leq x} \frac{1}{n^s} + \chi(s)\sum_{n \leq x} \frac{1}{n^{1-s}}
				+ (-1)^{[x]-1}e^{\pi i(s-1)/2} (2\pi t)^{s/2-1/2} e^{it/2-i\pi/8}	\nonumber\\
			 && \qquad \qquad \qquad \qquad  \times \Gamma(1-s) 
			\left( S_N + O \left( \left( \frac{CN}{t} \right)^{N/6} \right) + O(e^{-Ct}) \right),
			\label{RS_of_zeta}
\end{eqnarray}
where
\[
	S_N = \sum_{n=0}^{N-1} \sum_{\nu \leq n/2} \frac{n! i^{\nu - n}}{\nu ! (n-2\nu)! 2^n}
			\left(\frac{2}{\pi}\right)^{n/2-\nu} a_n \psi^{(n-2\nu)}\left( \frac{\eta}{\pi} - 2[x] \right),
\]
with $ a_n $ defined by
\[
	\exp{\left( (s-1) \log{\left( 1+\frac{z}{\sqrt{t}} \right)} -iz \sqrt{t} + \frac{1}{2}iz^2 \right)}
	= \sum_{n=0}^\infty a_n z^n,
\]
with $ a_0=1, a_n \ll t^{-n/2 + [n/3]} $.
R. Garunk\v{s}tis, A. Laurin\v{c}ikas, and J. Steuding proved an analogue of (\ref{RS_of_zeta})
for the Lerch zeta-function as follows;

\begin{theorem}[R. Garunk\v{s}tis, A. Laurin\v{c}ikas, and J. Steuding \cite{GLS}]\label{th:GLS}
Suppose that $ 0 < \alpha \leq 1, 0 < \lambda < 1 $ and $  0 \leq \sigma \leq 1 $. Suppose that 
$ t \geq 1, x = \sqrt{t/2\pi}, N = [x], M = [x - \alpha] $ and $ \beta = N - M $. Then
\begin{eqnarray}
	&& \zeta_L(s, \alpha, \lambda) 
		= \sum_{m = 0}^M \frac{e^{2\pi i m \lambda}}{(m + \alpha)^s} 
			+ \left( \frac{2\pi}{t} \right)^{\sigma -1/2 + it} e^{it+\pi i/4-2\pi i\{\lambda \} \alpha}
			\sum_{n=0}^N \frac{e^{-2\pi i \alpha n}}{(n+\lambda)^{1-s}}	\nonumber \\
	&& \qquad \qquad
		+ \left( \frac{2\pi}{t} \right)^{\sigma /2} e^{\pi i f(\lambda, \alpha, \sigma, t)}
		\phi(2x - 2N + \beta - \{ \lambda \} - \alpha) + O(t^{(\sigma-2)/2}),
\end{eqnarray}
where
\begin{eqnarray*}
	&& f(\lambda, \alpha, \sigma, t)=
	- \frac{t}{2\pi}\log{\frac{t}{2\pi e}} - \frac{7}{8} + \frac{1}{2} (\alpha^2 - \{\lambda\}^2)	\\
	&& \qquad \qquad \qquad \qquad \qquad 
		- \alpha \beta + 2y(\beta + \{\lambda\} - \alpha) - \frac{1}{2}(N+M) - \{ \lambda\}(\beta + \alpha).
\end{eqnarray*}
\end{theorem}

We prove an analogue of the approximate functional equation (\ref{AFE}) for (\ref{Hurwitz-zeta}) and 
(\ref{Lerch-zeta})
 (in Theorem \ref{th:Main_Theorem1}),
and gave another proof of the mean square formula for $ \zeta_L(1/2 + it, \alpha, \lambda) $ 
with respect to $ t $ (in Theorem \ref{th:Main_Theorem2}).

\begin{theorem}\label{th:Main_Theorem1}
Let $ 0 < \alpha \leq 1 $ and $ 0 < \lambda < 1 $.
Suppose that $ 0 \leq \sigma \leq 1,\ x \geq 1,\ y\geq 1 $ and $ 2\pi xy = |t| $. Then
\begin{eqnarray}
	&& \zeta_L(s, \alpha, \lambda) 
		= \sum_{0 \leq n \leq x} \frac{e^{2\pi i n \lambda}}{(n + \alpha)^s} \nonumber \\
	&& \qquad	+ \frac{\Gamma(1-s)}{(2\pi)^{1-s}}
		\left\{ e^{\{(1-s)/2 - 2 \alpha \lambda\}\pi i} 
				\sum_{0 \leq n \leq y}\frac{e^{2\pi in(1-\alpha)}}{(n+\lambda)^{1-s}} \right.	\nonumber \\
	&& \qquad \qquad \qquad \qquad \qquad \qquad \quad
		\left.	+e^{\{-(1-s)/2 + 2 \alpha (1-\lambda)\}\pi i } 
				\sum_{0 \leq n \leq y}\frac{e^{2\pi in \alpha}}{(n+1-\lambda)^{1-s}}
		\right\}	\nonumber \\
	&& \qquad + O(x^{-\sigma}) + O(|t|^{1/2-\sigma}y^{\sigma-1}). 		\label{AFE_of_Lerch}
\end{eqnarray}
Also, in the case $ \lambda = 1 $ that is $ \zeta_H(s, \alpha) $ it follows that 
\begin{eqnarray}
	&& \zeta_H(s, \alpha) = \sum_{0 \leq n \leq x}\frac{1}{(n+\alpha)^s} 	\nonumber \\
	&& \qquad \qquad \qquad \quad
		+ \frac{\Gamma(1-s)}{(2\pi)^{1-s}} 
		\left\{ e^{\frac{\pi i}{2}(1-s)} 
		\sum_{n \leq y} \frac{e^{2\pi in(1-\alpha)}}{n^{1-s}} + 
		e^{-\frac{\pi i}{2}(1-s)} \sum_{n \leq y} \frac{e^{2\pi in \alpha}}{n^{1-s}} \right\}	\nonumber \\
	&& \qquad \qquad + O(x^{-\sigma}) + O(|t|^{1-\sigma}y^{\sigma-1}).	\label{AFE_of_Hurwitz}
\end{eqnarray} 
\end{theorem}

\medskip

\begin{remark}
Theorem \ref{th:Main_Theorem1} can be proved by the method similar to the proof of
Theorem \ref{th:GLS}, but results of Theorem \ref{th:Main_Theorem1} has advantage of choosing 
parameters $ x $ and $ y $ freely, only under the condition $ 2\pi xy = |t| $ as compared with
the result of Theorem \ref{th:GLS}.
Also for approximate functional equations (\ref{AFE_of_Lerch}) and (\ref{AFE_of_Hurwitz}),
$ \zeta_L(s, \alpha, \lambda) $ is a generalization of $ \zeta_H(s, \alpha) $,
but (\ref{AFE_of_Lerch}) in Theorem \ref{th:Main_Theorem1}
does not include (\ref{AFE_of_Hurwitz}).
\end{remark}

\medskip

\begin{theorem}\label{th:Main_Theorem2}
Let $ 0 < \alpha \leq 1,\ 0 < \lambda \leq 1 $. Then,
\begin{equation}\label{mean_value_of_Lerch}
	\int_1^T \left|\zeta_L \left(\frac{1}{2} + it , \alpha, \lambda \right) \right|^2 dt
	= T \log{\frac{T}{2\pi}} + 
		\begin{cases}
			O(T(\log{T})^{1/2})	& (0 < \alpha <1 ),	\\
			O(T(\log{T})^{3/4})	& (\alpha = 1),
		\end{cases}	
\end{equation}
as $ T \rightarrow \infty $.
\end{theorem}

\medskip

\begin{remark}
The result of Theorem \ref{th:Main_Theorem2} has larger error term than the result already proved by
R. Garunk\v{s}tis, A. Laurin\v{c}ikas and J. Steuding \cite{GLS2}, and they proved using
Theorem \ref{th:GLS} (see \cite{GLS2}).
However, the main term on the right-hand side of (\ref{mean_value_of_Lerch}) can be obtained
more simply than the method of \cite{GLS2} by using Theorem \ref{th:Main_Theorem1}.
We will describe the proof of Theorem \ref{th:Main_Theorem2} in Section \ref{sec:proof_of_Thm3}.
\end{remark}

\section{Proof of Theorem \ref{th:Main_Theorem1}}

In this section, we prove Theorem \ref{th:Main_Theorem1}.
The basic tool of the proof is the same as the approximate functional equation for
the Riemann zeta-function (\ref{AFE}), that is the saddle point method.

\medskip

\textbf{Proof of Theorem \ref{th:Main_Theorem1}}. Let $ M \in \N $ be sufficiently large.
We have  
\begin{eqnarray}
	\zeta_L(s, \alpha, \lambda) 
	&=& \sum_{n=0}^M \frac{e^{2\pi i n \lambda}}{(n+\alpha)^s} 
						+ \sum_{n=M+1}^\infty \frac{e^{2\pi i n \lambda}}{(n+\alpha)^s} \nonumber\\
	&=& \sum_{n=0}^M \frac{e^{2\pi i n \lambda}}{(n+\alpha)^s} 
		+ \frac{e^{2\pi i \lambda M}}{\Gamma(s)}
		\int_0^\infty \frac{ t^{s-1} e^{(M+\alpha)t}}{e^{t-2\pi i \lambda}-1} dt \nonumber \\
	&=& \sum_{n=0}^M \frac{e^{2\pi i n \lambda}}{(n+\alpha)^s} 
		+ \frac{e^{2\pi i \lambda M}\Gamma(1-s)}{2\pi i e^{\pi i s}}
			\int_C \frac{z^{s-1}e^{(M+\alpha)z}}{e^{z-2\pi i \lambda}-1} dz,
		\label{contour of Hurwitz}
\end{eqnarray}
where $ C $ is the contour integral path that comes from $ +\infty $ to $ \varepsilon $ along 
the real axis, then goes along the circle of radius $ \varepsilon $ counter clockwise, and finally
goes from $ \varepsilon $ to $ +\infty $.

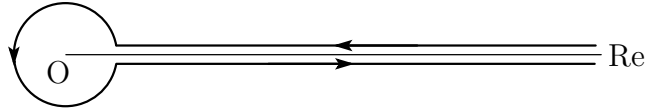
\begin{center}

{\unitlength 0.1in%
\begin{picture}( 30.9400,  5.4000)(  6.5000,-10.6000)%
%
\special{pn 13}%
\special{pa 2781 742}%
\special{pa 2361 742}%
\special{fp}%
\special{sh 1}%
\special{pa 2361 742}%
\special{pa 2428 762}%
\special{pa 2414 742}%
\special{pa 2428 722}%
\special{pa 2361 742}%
\special{fp}%
%
\special{pn 8}%
\special{pa 942 790}%
\special{pa 3744 790}%
\special{fp}%
%
\put(34.0000,-10.6000){\makebox(0,0)[lb]{}}%
\put(38.7500,-7.9500){\makebox(0,0){$ \mathrm{Re} $}}%
%
\special{pn 13}%
\special{pa 2002 838}%
\special{pa 2432 838}%
\special{fp}%
\special{sh 1}%
\special{pa 2432 838}%
\special{pa 2365 818}%
\special{pa 2379 838}%
\special{pa 2365 858}%
\special{pa 2432 838}%
\special{fp}%
\put(9.0000,-8.8000){\makebox(0,0){O}}%
%
\special{pn 13}%
\special{ar 940 790 270 270  0.1973956  6.1000745}%
%
\special{pn 13}%
\special{pa 670 778}%
\special{pa 670 828}%
\special{fp}%
\special{sh 1}%
\special{pa 670 828}%
\special{pa 690 761}%
\special{pa 670 775}%
\special{pa 650 761}%
\special{pa 670 828}%
\special{fp}%
\special{pa 670 828}%
\special{pa 670 828}%
\special{fp}%
%
\special{pn 13}%
\special{pa 1206 742}%
\special{pa 3715 742}%
\special{fp}%
%
\special{pn 13}%
\special{pa 1206 838}%
\special{pa 3710 838}%
\special{fp}%
\end{picture}}%

\end{center}

Let $ t>0 $ and $ x \leq y $, so that $ 1 \leq x \leq \sqrt{t/2\pi} $.
Let $ \sigma \leq 1, M = [x], N = [y], \eta = 2\pi y $. 
We deform the contour integral path $ C $ to the combination of the straight lines 
$ C_1, C_2, C_3, C_4 $ joining 
$ \infty,\ c\eta + i \eta(1+c) + 2\pi i \lambda,\ -c\eta + i\eta(1-c) + 2\pi i \lambda, 
  \ -c\eta-(2L+1)\pi i + 2\pi i \lambda ,\ \infty $,
where $ c $ is an absolute constant, $ 0 < c \leq 1/2 $.
 
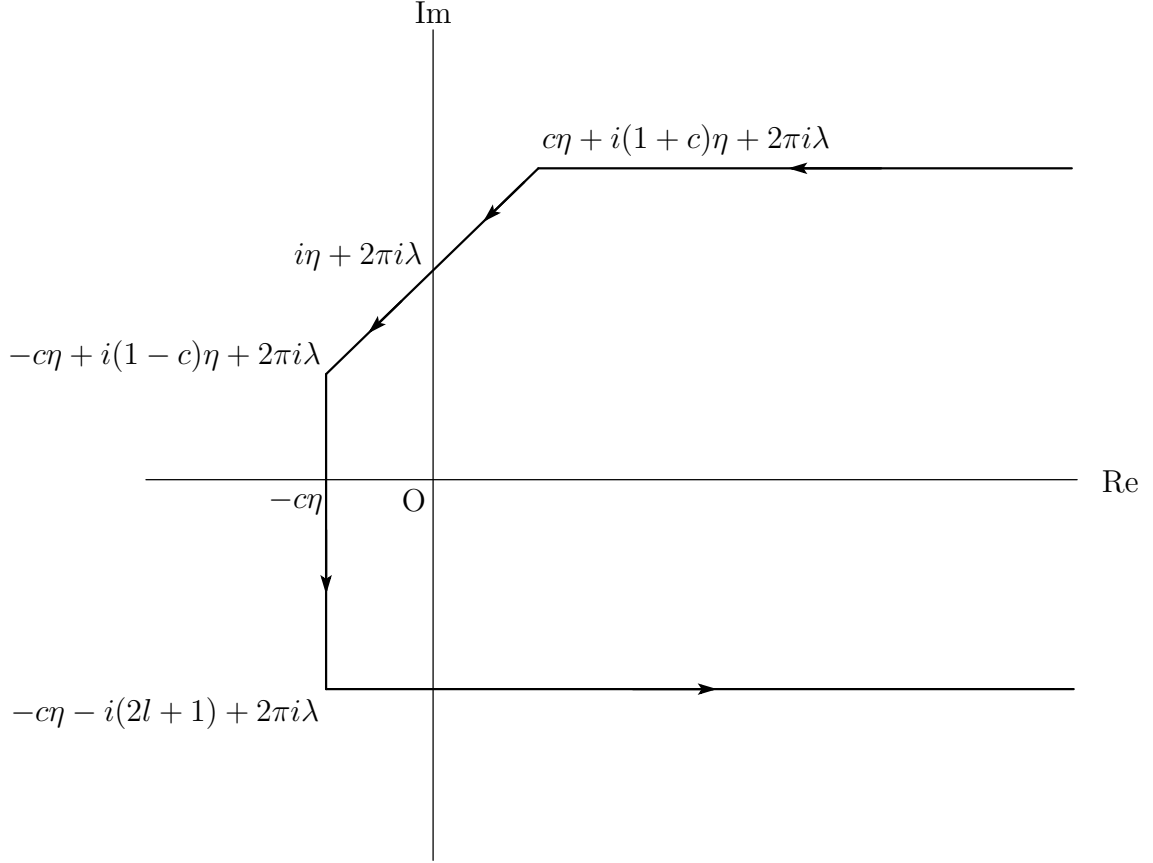
\begin{center}
{\unitlength 0.1in%
\begin{picture}( 59.0900, 45.0600)( -0.2900,-51.4000)%
%
\special{pn 8}%
\special{pa 2509 5140}%
\special{pa 2509 790}%
\special{fp}%
\special{pa 1006 3146}%
\special{pa 5880 3146}%
\special{fp}%
%
\special{pn 13}%
\special{pa 5862 4243}%
\special{pa 1949 4243}%
\special{fp}%
\special{pa 1949 4243}%
\special{pa 1949 2593}%
\special{fp}%
\special{pa 1949 2593}%
\special{pa 3060 1515}%
\special{fp}%
\special{pa 3060 1515}%
\special{pa 5852 1515}%
\special{fp}%
\put(19.2100,-25.8400){\makebox(0,0)[rb]{$ -c \eta + i(1-c)\eta + 2\pi i \lambda $}}%
\put(24.5300,-20.4900){\makebox(0,0)[rb]{$ i\eta + 2\pi i \lambda $}}%
\put(30.7900,-14.5100){\makebox(0,0)[lb]{$ c\eta + i(1+c)\eta + 2\pi i \lambda $}}%
\put(61.0500,-31.5500){\makebox(0,0){$ \rm{Re} $}}%
\put(19.1200,-42.8800){\makebox(0,0)[rt]{$ -c\eta -i(2l+1) + 2\pi i \lambda $}}%
\put(25.0900,-6.9900){\makebox(0,0){$ \rm{Im} $}}%
\put(24.0700,-32.5500){\makebox(0,0){O}}%
\put(17.9000,-32.6400){\makebox(0,0){$ -c\eta $}}%
%
\special{pn 13}%
\special{pa 4854 1515}%
\special{pa 4387 1515}%
\special{fp}%
\special{sh 1}%
\special{pa 4387 1515}%
\special{pa 4454 1535}%
\special{pa 4440 1515}%
\special{pa 4454 1495}%
\special{pa 4387 1515}%
\special{fp}%
%
\special{pn 13}%
\special{pa 2986 1587}%
\special{pa 2785 1782}%
\special{fp}%
\special{sh 1}%
\special{pa 2785 1782}%
\special{pa 2847 1750}%
\special{pa 2823 1745}%
\special{pa 2819 1721}%
\special{pa 2785 1782}%
\special{fp}%
%
\special{pn 13}%
\special{pa 2343 2211}%
\special{pa 2183 2366}%
\special{fp}%
\special{sh 1}%
\special{pa 2183 2366}%
\special{pa 2245 2334}%
\special{pa 2221 2329}%
\special{pa 2217 2305}%
\special{pa 2183 2366}%
\special{fp}%
%
\special{pn 13}%
\special{pa 3555 4243}%
\special{pa 3975 4243}%
\special{fp}%
\special{sh 1}%
\special{pa 3975 4243}%
\special{pa 3908 4223}%
\special{pa 3922 4243}%
\special{pa 3908 4263}%
\special{pa 3975 4243}%
\special{fp}%
%
\special{pn 13}%
\special{pa 1949 3409}%
\special{pa 1949 3726}%
\special{fp}%
\special{sh 1}%
\special{pa 1949 3726}%
\special{pa 1969 3659}%
\special{pa 1949 3673}%
\special{pa 1929 3659}%
\special{pa 1949 3726}%
\special{fp}%
\end{picture}}%
\end{center}
 
We calculate the residue of integrand of (\ref{contour of Hurwitz}). Since
\begin{eqnarray*}
	&& \lim_{z \rightarrow 2\pi i(\lambda + n)}\{z - 2\pi i(\lambda + n) \}
			\cdot \frac{z^{s-1} e^{-(M+\alpha)z}}{e^{z-2\pi i \lambda}-1}	\\
	&& = \lim_{z \rightarrow 2\pi i(\lambda + n)} 
			\left( \frac{e^{z - 2\pi i \lambda}-1}{z-2\pi i(\lambda + n)} \right)^{-1}
			e^{-(M+\alpha)z} \cdot z^{s-1}
	   = e^{- 2\pi i(M+\alpha)(\lambda + n)} (2\pi i(n + \lambda))^{s-1},
\end{eqnarray*}
we have
\begin{eqnarray*}
	&& \underset{z=2\pi i(\lambda + n)}{\rm{Res}} \frac{z^{s-1} e^{-(M+\alpha)z}}{e^{z-2\pi i \lambda}-1} 	\\
	&& \qquad \qquad =  e^{- 2\pi i(M+\alpha)(\lambda + n)} (2\pi (n + \lambda)i)^{s-1}		\\
	&& \qquad \qquad =
		\begin{cases}
			e^{- 2\pi i(M+\alpha)(\lambda + n)} (2\pi(n + \lambda) e^{\pi i/2})^{s-1} 	& (n \geq 0) 		\\
			e^{2\pi i(M+\alpha)(|n| - \lambda)} (2\pi(|n| - \lambda) e^{3\pi i/2})^{s-1}	& (n \leq -1) 
		\end{cases}		\\
	&& \qquad \qquad =
		\begin{cases}
			- \dfrac{e^{\pi i s}}{(2\pi)^{1-s}} \cdot e^{\{(1-s)/2 -2(M+\alpha) \lambda \}\pi i}
				\cdot \dfrac{e^{2\pi in(1-\alpha)}}{(n + \lambda)^{1-s}} 		& (n \geq 0) 		\\
			- \dfrac{e^{\pi i s}}{(2\pi)^{1-s}} \cdot e^{-\{(1-s)/2 +2(M+\alpha)(1-\lambda) \}\pi i}
				\cdot \dfrac{e^{2\pi i(-n)\alpha}}{(|n| - \lambda)^{1-s}}		& (n \leq -1) 
		\end{cases}
\end{eqnarray*}
and we have
\begin{eqnarray*}
	&& \sum_{n=-N+1}^N \underset{z=2\pi i n}{\rm{Res}}
		 \frac{z^{s-1} e^{-(M+\alpha)z}}{e^{z-2\pi i \lambda}-1} 	\\
	&& \quad = - \frac{e^{\pi is}}{(2\pi)^{s-1}}
		\left\{ e^{\{(1-s)/2 -2(M+\alpha) \lambda \}\pi i}
						\sum_{n=0}^N \frac{e^{2\pi in(1-\alpha)}}{(n + \lambda)^{1-s}} \right.	\\
	&& \qquad \qquad \qquad \qquad \qquad \qquad \left. + e^{-\{(1-s)/2 +2(M+\alpha)(1-\lambda) \}\pi i}
						\sum_{n=0}^N \frac{e^{2\pi in\alpha}}{(n+1-\lambda)^{1-s}} \right\}.
\end{eqnarray*}
Therefore we obtain
\begin{eqnarray}
	&&\zeta_L(s, \alpha, \lambda)
		= \sum_{n=0}^M \frac{e^{2\pi i n \lambda}}{(n+\alpha)^s}		\nonumber \\
	&&  \qquad + \frac{\Gamma(1-s)}{(2\pi)^{1-s}}
		\left\{ e^{\{(1-s)/2 -2\alpha \lambda \}\pi i}
						\sum_{n=0}^N \frac{e^{2\pi in(1-\alpha)}}{(n + \lambda)^{1-s}} \right.	\nonumber \\
	&& \qquad \qquad \qquad \qquad \qquad \qquad \qquad \left. + e^{-\{(1-s)/2 +2\alpha(1-\lambda) \}\pi i}
						\sum_{n=0}^N \frac{e^{2\pi in\alpha}}{(n+1-\lambda)^{1-s}} \right\}	\nonumber \\
	&& \qquad  
		+ \frac{e^{2\pi i \lambda M} \Gamma(1-s)}{2\pi i e^{\pi i s}}
			\left( \int_{C_1} + \int_{C_2} + \int_{C_3} + \int_{C_4} \right)
		\frac{z^{s-1}e^{-(M+\alpha)z}}{e^{z-2\pi i \lambda}-1}dz.
																				\label{zeta_L (1)}
\end{eqnarray}
From here, we consider the order of integral terms on right-hand side of (\ref{zeta_L (1)}).

First, we consider the integral path $ C_4 $.
Let $ z = u + iv = re^{i\theta } $ then $ |z^{s-1}| = r^{\sigma -1} $, and since 
$ \theta \geq  5\pi/4,\  r \gg \eta ,\  |e^{z-2\pi i \lambda} -1| \gg 1 $, we have
\begin{eqnarray}
	\int_{C_4} \frac{z^{s-1}e^{-(M+\alpha)z}}{e^{z - 2\pi i \lambda} - 1} dz
	& = & \int_{C_4} \frac{(re^{i\theta})^{\sigma+it-1}e^{-(M+\alpha )(u+iv)}}{e^{z - 2\pi i \lambda} - 1} dz	\nonumber\\
	& \ll & \eta^{\sigma-1} e^{-5\pi t/4} \int_{c\eta}^\infty e^{-(M+\alpha )u}du			\nonumber\\
	& = & \eta^{\sigma -1}(M + \alpha )^{-1}e^{(M+\alpha )c\eta - 5\pi t/4}					\nonumber\\
	& \ll & e^{(c-5\pi/4 )t}.		\label{evaluate C_4}
\end{eqnarray}

Secondly, we consider the order of integral on $ C_3 $ of (\ref{zeta_L (1)}).
Noting
\[
	\arctan{\varphi} = \int_0^\varphi \frac{d\mu}{1+\mu^2} > \int_0^\varphi \frac{d\mu}{(1+\mu)^2}
					 = \frac{\varphi}{1+\varphi}
\]
for $ \varphi > 0 $, we can write
\[
	\theta = \arg{z} = \frac{\pi}{2} + \arctan{\frac{c}{1-c}} = \frac{\pi}{2} + c + A(c)
\]
on $ C_3 $, where $ A(c) $ is a constant depending on $ c $. Then we have
\begin{eqnarray*}
	|z^{s-1}e^{-(M+\alpha)z}| 
	& = & r^{\sigma} e^{-t\theta + (M+\alpha)c\eta}		\\
	& \ll & \eta^{\sigma-1} e^{-(\pi/2 + c + A(c))t + (M+\alpha)\eta}	\\
	& \ll & \eta^{\sigma-1} e^{-(\pi/2 + A(c)}t.
\end{eqnarray*}
Therefore, since $ |e^{z-2\pi i \lambda}-1| \gg 1 $, we have
\begin{equation}
	\int_{C_3} \frac{z^{s-1} e^{-(M+\alpha )z}}{e^{z - 2\pi i \lambda} - 1} dz 
	\ll \eta^{\sigma} e^{-(\pi/2 + A(c))t}.		\label{evaluate C_3}	
\end{equation}

Thirdly, since $ |e^{z-2\pi i \lambda}-1| \gg e^u $ on $ C_1 $, we have
\[
	\frac{z^{s-1}e^{-(M+\alpha)z}}{e^{z - 2\pi i \lambda} - 1} 
	\ll \eta^{\sigma -1} \exp{\left(-t\arctan{\frac{(1+c)\eta}{u}} - (M+\alpha +1)u \right)}.
\]
Since $ M + \alpha + 1 \geq x = t/\eta $, the term $ (M + \alpha + 1)u $ on the 
right-hand side of the above may be replaced by $ tu/\eta $. Also, since
\[
	\frac{d}{du} \left( \arctan{\frac{(1+c)\eta}{u}} + \frac{u}{\eta} \right)
	= - \frac{(1+c)\eta}{u^2 + (1+c)^2 \eta^2} + \frac{1}{\eta} > 0
\]
we have
\[
	\arctan{\varphi} = \int_0^\varphi \frac{d\mu}{1+\mu^2} < \int_0^\varphi d\mu = \varphi.
\]
Therefore
\[
	\arctan{\frac{(1+c)\eta}{u}} + \frac{u}{\eta}
	\geq \arctan{\frac{1+c}{c}} + c = \frac{\pi}{2} - \arctan{\frac{c}{1+c}} + c
	> \frac{\pi}{2} + B(c)
\]
in $ u \geq c\eta $, where $ B(c) = c^2/(1+c)^2 $. Then we have
\[
	\frac{z^{s-1}e^{-(M+\alpha)z}}{e^{z - 2\pi i \lambda} - 1}
	\ll \eta^{\sigma -1} \exp{\left( - \left( \frac{\pi}{2} + B(c) \right)t \right)}.
\]
Since
\[
	\frac{z^{s-1}e^{-(M+\alpha)z}}{e^{z - 2\pi i \lambda} - 1}
	\ll
	\begin{cases}
		\eta^{\sigma -1} \exp{\left( - \left( \dfrac{\pi}{2} + B(c) \right)t \right)}  
				& (c\eta \leq u \leq \pi \eta), \\
		\eta^{\sigma -1} \exp{(-xu)} & (u \geq \pi \eta),
	\end{cases}
\]
we obtain
\begin{eqnarray}
	\int_{C_1} \frac{z^{s-1}e^{-(M+\alpha)z}}{e^{z - 2\pi i \lambda} - 1}
	&\ll & \eta^{\sigma -1} \left\{ \int_{c\eta}^{\pi \eta} 
		e^{- ( \pi/2 + B(c))t} du + \int_{\pi \eta}^\infty e^{-xu} du \right\}		\nonumber\\
	&\ll & \eta^\sigma e^{- ( \pi/2 + B(c))t} + \eta^{\sigma-1} e^{-\pi \eta x}		\nonumber\\
	&\ll & \eta^\sigma e^{- ( \pi/2 + B(c))t}	\label{evaluate C_1}.
\end{eqnarray}

Finally, we describe the evaluation of the integral on $ C_2 $.
Rewriting $ z = i(\eta + 2\pi \lambda ) + \xi e^{\pi i/4} $
(where $ \eta \in \R $ and $ |\eta| \leq \sqrt{2}c\eta $ ), we have 
\begin{eqnarray*}
	z^{s-1} &=& \exp\left\{ (s-1) \left( \log{(i(\eta + 2\pi \lambda) + \xi e^{-\pi i/4})} \right) \right\}	\\
			&=& \exp\left\{ (s-1) \left( \frac{\pi i}{2} 
				+ \log{(\eta + 2\pi \lambda + \xi e^{-\pi i/4})} \right) \right\}	\\
			&=& \exp\left\{ (s-1) \left( \frac{\pi i}{2} 
				+ \log({\eta}+ 2\pi \lambda) + \frac{\xi}{\eta + 2\pi \lambda}e^{-\pi i/4} \right. \right. \\
			& &	\left. \left. \qquad \qquad \qquad \qquad \qquad \qquad 
				- \frac{\xi^2}{2(\eta + 2\pi \lambda)^2}e^{-\pi i/2} + O\left( \frac{\xi^3}{\eta^3}\right) 
				\right) \right\}	\\
			&\ll& (\eta + 2\pi \lambda)^{\sigma-1} 
				\exp\left\{ \left( -\frac{\pi}{2} + \frac{\xi}{\sqrt{2}(\eta + 2\pi \lambda)} 
						- \frac{\xi^2}{2(\eta + 2\pi \lambda)^2} + O\left( \frac{\xi^3}{\eta^3}\right) 
				\right)t \right\}	
\end{eqnarray*}
as $ \eta \rightarrow \infty $. Also, since
\[
	 \frac{e^{-(M+\alpha)z}}{e^{z - 2\pi i \lambda} - 1}
	 = \frac{e^{-(M+\alpha-x)z}}{e^{z - 2\pi i \lambda} - 1} \cdot e^{-xz}
\]
and
\[
	\frac{e^{-(M+\alpha-x)z}}{e^{z - 2\pi i \lambda} - 1} \ll 
	\begin{cases}
		e^{(x-M-\alpha-1)u} & \left( u > \dfrac{\pi}{2} \right) \\
		e^{(x-M-\alpha)u}   & \left( u < -\dfrac{\pi}{2} \right),
	\end{cases}
\]	
we have
\[
	\frac{e^{-(M+\alpha)z}}{e^{z - 2\pi i \lambda} - 1}
	\ll |e^{-xz}| = e^{- \xi t/\sqrt{2} \eta}	\quad \left( |u| > \frac{\pi}{2} \right).
\]
Hence
\begin{eqnarray}
	&& \int_{C_2 \cap \{z \,|\, |u|>\pi/2\}} \frac{z^{s-1} e^{-(M+\alpha)z}}{e^{z - 2\pi i \lambda} - 1} dz	 	\nonumber\\
	&& \ll \int_{C_2 \cap \{z \,|\, |u|>\pi/2\}} 	
		(\eta + 2\pi \lambda)^{\sigma -1} \nonumber\\
	&& \quad 
		\times \exp\left\{ \left( -\frac{\pi}{2} + \frac{\xi}{\sqrt{2}(\eta + 2\pi \lambda)} 
						- \frac{\xi^2}{2(\eta + 2\pi \lambda)^2} + O\left( \frac{\xi^3}{\eta^3}\right) 
				\right)t \right\}
		\exp\left( - \frac{\xi t}{\sqrt{2} \eta} \right) d\xi	\nonumber\\
	&& \ll \int_{-\sqrt{2}c \eta}^{\sqrt{2}c \eta} 
		(\eta + 2\pi \lambda)^{\sigma -1} 
		e^{-\pi t/2} \exp\left\{ \left( - \frac{\xi^2}{2(\eta + 2\pi \lambda)^2} 
				+ O\left( \frac{\xi^3}{\eta^3}\right) \right)t \right\} d\xi	\nonumber\\
	&& \ll \int_{-\infty}^{\infty} 
		(\eta + 2\pi \lambda)^{\sigma -1} 
		e^{-\pi t/2} \exp\left\{ \left( - \frac{\xi^2}{2(\eta + 2\pi \lambda)^2} 
				+ O\left( \frac{\xi^3}{\eta^3}\right) \right)t \right\} d\xi	\nonumber\\
	&& \ll \eta^{\sigma -1} e^{- \pi t/2} \int_{-\infty}^\infty 
		\exp\left\{ - \frac{D(c)\xi^2 t}{\eta^2} \right\} d\xi	\nonumber\\
	&& \ll \eta^\sigma t^{-1/2} e^{\pi t/2}, \label{evaluate C_2_1}
\end{eqnarray}
where $ D(c) $ is a constant depending on $ c $.
The argument can also be applied to the part $ |u| \leq \pi/2 $ if $ |e^{z-2\pi i \lambda}| > A $. 
If not, that is the case when the contour goes too near to the pole at $ z = 2 \pi i N + 2\pi i \lambda $,
we take an arc of the circle $ |z - 2 \pi i N - 2\pi i \lambda | = \pi/2 $.
On this arc we can write to $ z = 2\pi i N + 2\pi \lambda i + (\pi/2)e^{i \beta} $, and
\begin{eqnarray*}
	\log{(z^{s-1})} 
	&=& (s-1) \log{\left( 2\pi i N + 2\pi i \lambda + \frac{\pi}{2}e^{i \beta} \right)}	\\
	&=& (s-1) \log{e^{\pi i/2} \left( 2\pi N + 2\pi \lambda 
							+ \frac{\pi}{2} \cdot \frac{e^{i \beta}}{i} \right)}	\\
	&=& (\sigma + it - 1) \left\{ \frac{\pi i}{2} + \log(2\pi(N+\lambda)) 
							+ \log\left( 1 + \frac{e^{i \beta}}{4(N+\lambda)i} \right) \right\}	\\
	&=& -\frac{\pi t}{2} + (s-1)\log(2\pi(N+\lambda)) + \frac{te^{i\beta}}{4(N+\lambda)} + O(1).
\end{eqnarray*}
On the last line of the above calculations, we used $ N^2 \gg t $ which follows from the assumption $ x \leq y $.
Then
\begin{eqnarray*}
	&& z^{s-1} e^{-(M+\alpha)z} \\
	&& \qquad 
		= \exp{ \left( -\frac{\pi}{2} + (s-1)\log(2\pi(N+\lambda))+ \frac{te^{i\beta}}{4(N+\lambda )}
		- \frac{\pi}{2}(M+\alpha )e^{i \beta} + O(1)\right)},
\end{eqnarray*}
and since
\[
	\frac{te^{i\beta}}{4(N+\lambda)} - \frac{\pi}{2}(M+ \alpha)e^{i\beta}
	= \frac{2\pi xy - 2 \pi([x] + \alpha)([y] + \lambda)}{4(N + \lambda)} e^{i \beta} = O(1)
\]
we have
\begin{eqnarray*}
	z^{s-1} e^{-(M+\alpha)z}
	& \ll & \exp{ \left( -\frac{\pi}{2} + (s-1)\log(2\pi(N+\lambda)) + O(1) \right)}	\\
	& \ll & N^{\sigma - 1} e^{-\pi t/2}.
\end{eqnarray*}
Hence, the integral on the small semicircle can be evaluated as 
$ O(\eta^{\sigma - 1} e^{-\pi t/2}) $. Therefore together with (\ref{evaluate C_2_1}),
we have
\begin{equation}
	\int_{C_2} \frac{z^{s-1} e^{-(M+\alpha)z}}{e^{z - 2\pi i \lambda} - 1} dz
	\ll \eta^\sigma t^{-1/2}e^{-\pi t/2} + \eta^{\sigma -1} e^{-\pi t/2}.
	\label{evaluate C_2}
\end{equation}

Now, evaluation of all the integrals was done. Using the results (\ref{evaluate C_4}), 
(\ref{evaluate C_3}), (\ref{evaluate C_1}), (\ref{evaluate C_2}) and 
$ e^{2\pi i (\lambda N - s/2)} \Gamma(1-s) \ll t^{1/2 - \sigma} e^{\pi t/2} $, 
we see that the integral term of (\ref{zeta_L (1)}) is
\begin{eqnarray*}
	&\ll& t^{1/2 - \sigma} e^{\pi t/2} 
		\{ \eta^\sigma e^{-(\pi/2 + B(c))t}
			+ \eta^\sigma t^{-1/2} e^{-\pi t/2} + \eta^{\sigma-1}e^{\pi t/2} 	\\	
	&& \qquad \qquad \qquad \qquad \qquad \qquad \qquad \qquad
			+ \eta^\sigma e^{t(\pi/2 + A(c))} + e^{(c+ \varepsilon - 5\pi/4)t} \}  \\
	&\ll& t^{1/2} \left( \frac{\eta}{t}\right)^\sigma e^{-(A(c)+B(c))t} 
			+ \left(\frac{\eta}{t}\right)^\sigma + t^{-1/2} \left(\frac{\eta}{t}\right)^{\sigma-1}
			+ t^{1/2 - \sigma} e^{(c+ \varepsilon - 5\pi/4)t}	\\
	&\ll& e^{-\delta t} + x^{-\sigma} + t^{-1/2} x^{1-\sigma}  \ll  x^{-\sigma},
\end{eqnarray*}
where $ \delta $ is a small positive real number.
Therefore we have
\begin{eqnarray}
	&& \zeta_L(s, \alpha, \lambda) 
		= \sum_{0 \leq n \leq x} \frac{e^{2\pi i n \lambda}}{(n + \alpha)^s} \nonumber \\
	&& \qquad	+ \frac{\Gamma(1-s)}{(2\pi)^{1-s}}
		\left\{ e^{\{(1-s)/2 - 2 \alpha \lambda\}\pi i} 
				\sum_{0 \leq n \leq y}\frac{e^{2\pi in(1-\alpha)}}{(n+\lambda)^{1-s}} \right.	\nonumber \\
	&& \qquad \qquad \qquad \qquad \qquad \qquad
		\left.	+e^{\{-(1-s)/2 + 2\alpha (1-\lambda)\}\pi i } 
				\sum_{0 \leq n \leq y}\frac{e^{2\pi in \alpha}}{(n+1-\lambda)^{1-s}}
		\right\}	\nonumber \\
	&& \qquad  + O(x^{-\sigma}),	\label{x_leq_y}
\end{eqnarray} 
that is, Theorem \ref{th:Main_Theorem1} in the case of $ x \leq y $ has been proved.

To prove Theorem \ref{th:Main_Theorem1} in the case $ x \geq y $, we use 
the following functional equation of the Lerch zeta-function;
\begin{eqnarray}
	&& \zeta_L(s, \alpha, \lambda) = \frac{\Gamma(1-s)}{(2\pi)^{1-s}}
		\{ e^{\{(1-s)/2 - 2 \alpha \lambda\}\pi i} 
				\zeta_L(1-s, \lambda, 1-\alpha) 	\nonumber \\
	&& \qquad \qquad \qquad \qquad \qquad
		 + e^{\{-(1-s)/2 + 2\alpha (1-\lambda)\}\pi i } 
				\zeta_L(1-s, 1-\lambda, \alpha) \}.
				\label{FE_Lerch}
\end{eqnarray}
Applying (\ref{x_leq_y}) to $ \zeta_L(1-s, \lambda, 1-\alpha) $ and $ \zeta_L(1-s, 1-\lambda, \alpha) $,
and substitute these into (\ref{FE_Lerch}), we have
\begin{eqnarray}
	&& \zeta_L(s, \alpha, \lambda) \nonumber \\
	&&= \frac{\Gamma(1-s)}{(2\pi)^{1-s}}
		\left[ e^{\{(1-s)/2 - 2 \alpha \lambda\}\pi i} 
			\left\{\sum_{0 \leq n \leq x} \frac{e^{2\pi i n \lambda}}{(n + \lambda)^{1-s}} \right.\right.	
			\nonumber \\
	&& \quad \left.
		+ \frac{\Gamma(s)}{(2\pi)^s}
		\left( e^{\{s/2 - 2\lambda(1- \lambda) \}\pi i} 
				 \sum_{0 \leq n \leq y}\frac{e^{2\pi in(1-\lambda)}}{(n+1-\alpha)^s} 
		+ e^{\{-s/2 + 2\alpha \lambda \}\pi i } 
		\sum_{0 \leq n \leq y} \frac{e^{2\pi i n \lambda}}{(n + \alpha)^s} \right)\right\}	\nonumber \\
	&& \qquad \qquad \qquad
		 + e^{\{-(1-s)/2 + 2\alpha (1-\lambda)\}\pi i }  
				 \left\{ \sum_{0 \leq n \leq x} \frac{e^{2\pi i n \alpha}}{(n + 1- \lambda)^{1-s}} \right.  
				\nonumber \\
	&& \quad	+ \frac{\Gamma(s)}{(2\pi)^s}
		\left( e^{\{(s/2 - 2 (1-\lambda) \alpha \}\pi i} 
				\sum_{0 \leq n \leq y}\frac{e^{2\pi in \lambda}}{(n+\alpha)^s} \right.	 \nonumber \\
	&& \qquad \qquad \qquad \qquad \qquad \qquad \qquad
		\left.\left.\left.	+e^{\{-s/2 + 2(1-\lambda)(1-\alpha)\}\pi i } 
				\sum_{0 \leq n \leq y}\frac{e^{2\pi in (1-\lambda)}}{(n+\alpha)^s}
		\right)\right\} \right]	\nonumber\\
	&& \quad + O(\Gamma(\sigma -1)(2\pi)^{-\sigma}(e^{\pi t/2}+e^{-\pi t/2})x^{\sigma -1})	\nonumber \\
	&& = \sum_{0 \leq n \leq y} \frac{e^{2\pi i n \lambda}}{(n + \alpha)^s} 
		+ \frac{\Gamma(1-s)}{(2\pi)^{1-s}}
		\left\{ e^{\{(1-s)/2 - 2 \alpha \lambda\}\pi i} 
				\sum_{0 \leq n \leq x}\frac{e^{2\pi in(1-\alpha)}}{(n+\lambda)^{1-s}} \right.	\nonumber \\
	&& \qquad \qquad \qquad \qquad \qquad \qquad 
		\left.	+e^{\{-(1-s)/2 + 2\alpha (1-\lambda)\}\pi i } 
				\sum_{0 \leq n \leq x}\frac{e^{2\pi in \alpha}}{(n+1-\lambda)^{1-s}}
		\right\}  \nonumber \\
	&& \quad + O(t^{1/2-\sigma}x^{\sigma-1}).	\nonumber
\end{eqnarray}
Interchanging $ x $ and $ y $, we obtain the theorem with $ x \geq y $. 
Combining this equation with (\ref{x_leq_y}), we obtain the proof of (\ref{AFE_of_Lerch}).

The proof of (\ref{AFE_of_Hurwitz}) is similar. 
However, the four integral path $ C_1, C_2, C_3 $ and $ C_4 $ are different from the proof of
(\ref{AFE_of_Lerch}), that is, as follows; 
The straight lines $ C_1, C_2, C_3, C_4 $ joining 
$ \infty,\ c\eta + i \eta(1+c),\ -c\eta + i\eta(1-c),\ -c\eta-(2L+1)\pi i,\ \infty $,
where $ c $ is an absolute constant, $ 0 < c \leq 1/2 $. 
Also, in the proof for the case $ x \geq y $, we use the functional equation
\begin{eqnarray*}
	 \zeta_H(s, \alpha) = \frac{\Gamma(1-s)}{(2\pi)^{1-s}}
		\{ e^{(1-s)\pi i/2} \zeta_L (1-s, 1, 1-\alpha ) 
			+ e^{-(1-s)\pi i/2} \zeta_L (1-s, 1, \alpha ) \}  \label{FE_Hurwitz},
\end{eqnarray*}
but this equation is not included in the functional equation (\ref{FE_Lerch}).
Noticing these points, we can prove (\ref{AFE_of_Hurwitz}) by a similar method.
This completes the proof of Theorem \ref{th:Main_Theorem1}.
\qed

\section{Proof of Theorem \ref{th:Main_Theorem2}} \label{sec:proof_of_Thm3}
In this section, using Theorem \ref{th:Main_Theorem1}, we give the proof of 
Theorem \ref{th:Main_Theorem2}.

\medskip

\textbf{Proof of Theorem \ref{th:Main_Theorem2}}. 
Let
\[
	x = \frac{t}{2\pi \sqrt{\log{t}}}, \quad y = \sqrt{\log{t}}
\]
and we assume $ t > 0 $ satisfies $ x \geq 1 $ and $ y \geq 1 $.
Use the Stirling formula
\[
	\Gamma(1-s)	e^{\{(1-s)/2 -2\alpha \lambda \}\pi i} \ll 1, \ 
	\Gamma(1-s) e^{\{-(1-s)/2 -2\alpha (1-\lambda) \}\pi i} \ll 1.
\]
Then if $ 0 < \lambda < 1 $, using (\ref{AFE_of_Lerch}) we have
\begin{eqnarray}
	&& \zeta_L \left(\frac{1}{2} + it, \alpha, \lambda \right) 
		= \sum_{0 \leq n \leq x} \frac{e^{2\pi i n \lambda}}{(n + \alpha)^{1/2 + it}} \nonumber	\\
	&& \qquad \qquad \qquad \qquad \qquad 
			+ O \left(\sum_{0 \leq n \leq y}\frac{e^{2\pi in(1-\alpha)}}{(n+\lambda)^{1/2 - it}} 
			+ \sum_{0 \leq n \leq y}\frac{e^{2\pi in \alpha}}{(n+1-\lambda)^{1/2 - it}} \right)	\nonumber \\
	&& \qquad \qquad \qquad \qquad \qquad 
			+ O(t^{-1/2}(\log{t})^{1/4}) + O((\log{t})^{-1/4}), 	\label{Order_of_Lerch}
\end{eqnarray}
and if $ \lambda = 1 $, using (\ref{AFE_of_Hurwitz}) we have
\begin{eqnarray}
	&& \zeta_H \left(\frac{1}{2} + it, \alpha \right) 
		= \sum_{0 \leq n \leq x} \frac{1}{(n + \alpha)^{1/2 + it}} 
			+ O \left(\sum_{n \leq y}\frac{e^{2\pi in(1-\alpha)}}{n^{1/2 - it}} 
			+ \sum_{n \leq y}\frac{e^{2\pi in \alpha}}{n^{1/2 - it}} \right)	\nonumber \\
	&& \qquad \qquad \qquad \qquad \qquad 
			+ O(t^{-1/2}(\log{t})^{1/4}) + O((\log{t})^{-1/4}) 	\label{Order_of_Hurwitz}.
\end{eqnarray}
\begin{enumerate}
	\item[(i)] In the case $ 0 < \lambda < 1 $ and $ 0 < \alpha < 1 $, since
	\[
		\sum_{n=0}^\infty \frac{e^{2\pi in(1-\alpha)}}{(n+\lambda)^{1/2}}, \quad 
			 \sum_{n=0}^\infty\frac{e^{2\pi in \alpha}}{(n+1-\lambda)^{1/2}}
	\]
	are convergent, and $ t^{-1/2}(\log{t})^{1/4} = o(1) $, $ (\log{t})^{-1/4} = o(1) $, 
	we have
	\[
		\zeta_L \left(\frac{1}{2} + it, \alpha, \lambda \right) 
		= \sum_{0 \leq n \leq x} \frac{e^{2\pi i n \lambda}}{(n + \alpha)^{1/2 + it}} + O(1).
	\]
	\item[(ii)] In the case $ 0 < \lambda < 1 $ and $ \alpha = 1 $, the second term on 
	right-hand side of (\ref{AFE_of_Lerch}) is
	\[
		\ll \int_0^y \frac{1}{(u+\lambda)^{1/2}}du = O(\sqrt{y}) = O((\log{t})^{1/4}),
	\]
	so we have
	\[
		\zeta_L \left(\frac{1}{2} + it, 1, \lambda \right) 
		= \sum_{n \leq x} \frac{e^{2\pi i n \lambda}}{n^{1/2 + it}} 
			+ O((\log{t})^{1/4}).
	\]
	\item[(iii)] In the case $ \lambda = 1 $ and $ 0 < \alpha < 1 $,
	consider similarly as in the case of (i) to obtain
	\[
		\zeta_L \left(\frac{1}{2} + it, \alpha , 1 \right) 
		= \sum_{0 \leq n \leq x} \frac{1}{(n+ \alpha)^{1/2 + it}} + O(1).
	\]
	\item[(iv)] In the case $ \lambda = 1 $ and $ \alpha = 1 $, since
	$ \zeta_L(s, 1, 1) = \zeta(s) $ we obtain
	\[
		\zeta_L \left(\frac{1}{2} + it, 1, 1\right) 
		= \sum_{n \leq x} \frac{1}{n^{1/2 + it}} + O((\log{t})^{1/4})
	\]
	(see Chap. VII in \cite{Tit}).
\end{enumerate}
Let
\[
	\Sigma(\alpha, \lambda) = \sum_{0 \leq n \leq x }\frac{e^{2\pi i n \lambda}}{(n + \alpha)^{1/2 + it}},
\]
and calculate as
\begin{eqnarray*}
	|\Sigma(\alpha, \lambda)|^2 
	&=& \mathop{\sum\sum}\limits_{0 \leq m, n \leq x}
				\frac{e^{2\pi i (m-n) \lambda}}{(m+ \alpha)^{1/2}(n+\alpha)^{1/2}} 
				\left( \frac{n+\alpha}{m+\alpha} \right)^{it}		\\
	&=& \sum_{0 \leq n \leq x}\frac{1}{n + \alpha} + 
		\mathop{\mathop{\sum\sum}\limits_{0 \leq m, n \leq x}}\limits_{m \neq n}
		\frac{e^{2\pi i (m-n) \lambda}}{(m+ \alpha)^{1/2}(n+\alpha)^{1/2}}
		\left( \frac{n+\alpha}{m+\alpha} \right)^{it}.
\end{eqnarray*}
Also $ T_1 = T_1(m,n) $ is a function in $ m, n $ satisfying
\[
	\max\{m,n\} = \frac{T_1}{2\pi \sqrt{\log{T_1}}}.
\]
Let $ X = T/2\pi \sqrt{\log{T}} $, then 
\begin{eqnarray}
	&& \int_1^T |\Sigma(\alpha, \lambda)|^2 dt
		= \sum_{0 \leq n \leq X}\frac{1}{n + \alpha} \{ T - T_1(n,n) \}	 \nonumber	\\
	&& \qquad \qquad \qquad \quad
		+ O \left( \mathop{\sum\sum}\limits_{0 \leq m < n \leq X}
				\frac{e^{2\pi i (m-n) \lambda}}
					{(m+ \alpha)^{1/2}(n+\alpha)^{1/2}} 
					\left( \log \frac{n+\alpha}{m+\alpha} \right)^{-1}	\right).
																		\label{int_Sigma^2}
\end{eqnarray}
Here, since
\[
	n \sqrt{\log{n}} 
	= \frac{T_1}{2\pi \sqrt{\log{T_1}}} \left( \log\frac{T_1}{2\pi \sqrt{\log{T_1}}} \right)^{1/2}
			\sim \frac{1}{2\pi} T_1(n,n)
\]
and
\begin{eqnarray*}
	\mathop{\sum\sum}\limits_{0 \leq m < n \leq X}
				\frac{e^{2\pi i (m-n) \lambda}}
					{(m+ \alpha)^{1/2}(n+\alpha)^{1/2}} 
					\left( \log \frac{n+\alpha}{m+\alpha} \right)^{-1}	
	\ll  X \log{X}	\ll T (\log{T})^{1/2}
\end{eqnarray*}
(see Lemma 3 in \cite{GLS2} or Lemma 2.6 in \cite{LG}), (\ref{int_Sigma^2}) can be rewritten as
\begin{equation}
	\int_1^T |\Sigma(\alpha, \lambda)|^2 dt
	= T \log{\frac{T}{2\pi}} + O(T(\log{T})^{1/2}). 		\label{int_Sigma^2_1}
\end{equation}
Therefore from (i), (ii), (iii), (iv) and (\ref{int_Sigma^2_1})
, and the Cauchy-Schwarz inequality, we obtain
\begin{eqnarray*}
	&& \int_1^T \left| \zeta_L \left(\frac{1}{2} + it, \alpha, \lambda \right) \right|^2 dt  \\
	&& = \int_1^T |\Sigma(\alpha, \lambda)|^2 dt +
		\begin{cases}
			O(T^{1/2} (\log{T})^{1/4})) + O(T) 			& (0 < \alpha <1 ),	 	\\
			O(T (\log{T})^{3/4}) + O(T (\log{T})^{1/2})	& (\alpha = 1)
		\end{cases}	
		\nonumber \\
	&& = T \log{\frac{T}{2\pi}} + 
		\begin{cases}
			O(T(\log{T})^{1/2})	& (0 < \alpha <1 ),	\\
			O(T(\log{T})^{3/4})	& (\alpha = 1).
		\end{cases}	
\end{eqnarray*}
Thus we obtain the proof of Theorem \ref{th:Main_Theorem2}.
\qed

\bigskip \bigskip
\author{Takashi Miyagawa}:	\\
Graduate School of Mathematics, \\
Nagoya University, \\
Chikusa-ku, Nagoya, 464-8602 Japan	\\
E-mail: d15001n@math.nagoya-u.ac.jp

\bigskip

\begin{thebibliography}{99}
\bibitem{GLS}	R. Garunk\v{s}tis, A. Laurin\v{c}ikas, and J. Steuding,
				An approximate functional equation for the Lerch zeta-function, 
				Mathematical Notes, Vol. \textbf{74}, No. 4 (2003), pp. 469-476.
\bibitem{GLS2}	R. Garunk\v{s}tis, A. Laurin\v{c}ikas, and J. Steuding,
				On the mean square of Lerch zeta-function, Arch. Math. (Basel) \textbf{80} (2003),
				no. 1, 47-60.
\bibitem{LG0}	A. Laurin\v{c}ikas, R. Garunk\v{s}tis, 
				On the Lerch zeta-function, Lith. Math. J. \textbf{36} (1996), 337-346.
\bibitem{LG}	A. Laurin\v{c}ikas, R. Garunk\v{s}tis, The Lerch zeta-function, Kluwer Academic Publishers,
				Dordrecht, 2002.
\bibitem{Tit}	E. C. Titchmarsh, The Theory of the Riemann Zeta-function. 2nd ed., Edited and
				with a preface by D. R. Heath-Brown, The Clarendon Press, Oxford University
				Press, New York, 1986.

\end{thebibliography}
\end{document}